\documentclass[11pt]{amsart}
\usepackage[margin=1.12in]{geometry}
\usepackage{amsmath,amssymb,amsthm,mathtools,booktabs,array,longtable}
\usepackage{microtype}
\usepackage{xurl}
\usepackage[hidelinks]{hyperref}
\hypersetup{pdftitle={Complete self-shrinkers with scalar curvature bounded below}}

\numberwithin{equation}{section}
\newtheorem{theorem}{Theorem}[section]
\newtheorem{proposition}{Proposition}[section]
\newtheorem{conjecture}{Conjecture}[section]
\newtheorem{lemma}[theorem]{Lemma}

\theoremstyle{remark}
\newcommand{\R}{\mathbb R}
\newcommand{\Lop}{\mathcal L}
\newcommand{\tr}{\operatorname{tr}}
\allowdisplaybreaks[2]
\title[Complete self-shrinkers ]{Complete self-shrinkers with \\ scalar curvature bounded below}
\author [Q. -M. Cheng, F. Li  and G. Wei]{Qing-Ming Cheng, Fengjiang Li{$^+$}  and Guoxin Wei{$^+$}\\  \\
{\smaller Dedicated to Professor Yuanlong Xin  for his 85th birthday}}

\address{Qing-Ming Cheng \newline
\indent
Mathematical Science Research Center, \newline 
\indent Chongqing University of Technology, Chongqing 400054, P. R. China. 
\vskip 1mm
\indent
Osaka Central  Advanced Mathematical Institute \newline
\indent
Osaka Metropolitan University, Osaka 558-8585, Japan. \newline
\indent qingmingcheng@yahoo.com, chengqingming@cqut.edu.cn }
\address{Fengjiang Li \newline
\indent
Mathematical Science Research Center, \newline 
\indent Chongqing University of Technology, Chongqing 400054, P. R. China. \newline 
\indent fengjiangli@cqut.edu.cn }

\address{Guoxin Wei \newline
\indent   School of Mathematical Sciences, \newline 
\indent  South China Normal University, Guangzhou 510631, P. R. China. \newline 
\indent weiguoxin@tsinghua.org.cn}

\thanks{$+$ Fengjiang Li  and Guoxin Wei are corresponding authors}
\thanks{*This work was supported by  National Natural Science Foundation of China (Grant Nos.
12301062, 12571050, 12671065), Natural Science Foundation of Chongqing
(No.CSTB2024NSCQ-MSX0537), Japan Society for the Promotion of Science Grant-in-Aid
for Scientific Research (C) (Grant No. 25K06992).}
\date{}

\begin{document}
\begin{abstract}
Let $X:M^n\to\mathbb R^{n+1}$ be a complete self-shrinker with constant
squared norm of the second fundamental form. We prove that 
up to an orthogonal transformation, the image $X(M^n)$
is  a sphere  $S^n(\sqrt n)$, or a hyperplane $\mathbb R^n$ or 
a generalized cylinder
$ S^k(\sqrt{k})\times\R^{n-k},\  1\le k\le n-1$ if scalar curvature $R$ satisfies  $R\geq -8/5$, which improves the recent results in \cite{CLW1}. 
\end{abstract}
\maketitle

\section{Introduction}
Self-shrinkers arise as models for singularities of mean curvature flow. We consider a connected 
complete immersed hypersurface \(X:M^n\to\R^{n+1}\), without boundary, satisfying
\begin{equation}\label{shrinker}
H+\langle X,N\rangle=0,
\end{equation}
where \(H\) is the mean curvature and \(N\) is a unit normal vector. We write \(\vec h\) for the second 
fundamental form, \(\lambda_i\) for its principal curvatures, and
\[
 S=|\vec h|^2,\qquad f_k=\sum_i\lambda_i^k.
\]
Self-shrinkers are closely related to singularities of the mean curvature
flow; see Huisken \cite{Huisken} and Colding--Minicozzi \cite{CM}.
In this paper we consider complete immersed self-shrinkers for which $S$
is constant. No assumption on polynomial volume growth or embeddedness
is imposed.

For a hyperplane through the origin, $S=0$. For the sphere
$S^n(\sqrt n)$ and the cylinders $S^k(\sqrt k)\times\mathbb R^{n-k}$,
$S=1$. It is natural to ask whether these are the only complete
self-shrinkers with constant $S$. This question is related to the gap
problem for minimal hypersurfaces in spheres; see Simons \cite{s},
Peng--Terng \cite{pt} and Yang--Cheng \cite{yc}.
We recall the following conjecture.

\begin{conjecture}
For an $n$-dimensional  complete self-shrinker \(X:M^n\to\R^{n+1}\) with constant squared norm \(S\) of the second fundamental form,
up to an orthogonal transformation, the image $X(M^n)$ is  a hyperplane $\mathbb R^n$ through the origin, a round sphere \(S^n(\sqrt n)\), or a generalized cylinder
\[
 S^k(\sqrt{k})\times\R^{n-k},\qquad 1\le k\le n-1.
\]
\end{conjecture}

The conjecture is resolved  in dimension two by Cheng and Ogata \cite{CO}. In dimension three, Cheng, Li and Wei \cite{CLWf3, CLWf4}
proved it when \(f_3\) or \(f_4\) is also constant; Recently, Fu and Yang \cite{FY}  and Guan and Zhao \cite{GZ}  have announced that they 
remove the condition  that \(f_3\) or \(f_4\) is  constant for dimension three. 
For further related research on pinching results, see  Cao and Li \cite{CaoLi}, 
Cheng and Peng \cite{CP},  Cheng and Wei \cite{CW} , 
Ding and Xin \cite{DX}, Lei, Xu and Xu \cite{LeiXuXu}, and Cheng, Wei and Yano \cite{CWY}. 
The relation between polynomial  volume growth and properness was studied  by Ding and Xin \cite{DXVolume} 
and Cheng and Zhou \cite{ChengZhou}. 
The generalized maximum principle of 
Cheng and Peng \cite{CP} permits a different approach without an assumption of polynomial volume growth. 
Recently,  Cheng, Li and Wei \cite{CLW1} have resolved the above conjecture if scalar curvature $R$ satisfies $R\geq -\frac75$.
For work under constant scalar curvature assumptions, see also \cite{CLW1, CLW2, g}.

In this paper,
our purpose is to study the above conjecture. We obtain the following:
\begin{theorem}\label{main}
For an $n$-dimensional  complete self-shrinker \(X:M^n\to\R^{n+1}\) with constant squared norm \(S\) of the second fundamental form,
up to an orthogonal transformation, the image $X(M^n)$ is  a hyperplane $\mathbb R^n$ through the origin, a round sphere \(S^n(\sqrt n)\), 
or a generalized cylinder
\[
 S^k(\sqrt{k})\times\R^{n-k},\qquad 1\le k\le n-1.
\]
if scalar curvature $R$ satisfies $R\geq -\frac85$.
\end{theorem}
The proof uses the fixed auxiliary function:
\[
F=\frac S4f_4-\frac{137}{600}f_3^2+\frac{f_3^4}{120S^3}.
\]
In contrast to minimal hypersurfaces in a sphere, the mean curvature of a self-shrinker need not vanish. We use the projection in $\{\vec e_1,\vec\psi\}^{\perp}$ and derive the necessary estimates without the trace-free condition. The condition $|\nabla F|\to0$ in the generalized maximum principle plays an important role in our proof.

\vskip5mm
\noindent
{\bf Organization}. 
Section 2  records the geometric identities and fixes notation.
Section 3  proves  key lemmas.  Section 4 derives main estimates
and Section 5  proves Theorem~\ref{main}.

\section{Preliminaries}
We choose an orthonormal frame at each point so that $h_{ij}=\lambda_i\delta_{ij}$, where 
\(\lambda_i\) are  principal curvatures, and set
$$
\begin{aligned}
f_k&=\sum_i\lambda_i^k,  \quad f=f_4-\frac{f_3^2}{S}, \quad A=\sum_{ijk}\lambda_i^2h_{ijk}^2,\\
B&=\sum_{ijk}\lambda_i\lambda_jh_{ijk}^2, \quad C =\sum_{ijk}\lambda_i h_{ijk}^2,
\quad D_k=\sum_i\lambda_i^2h_{iik}.
\end{aligned}
$$
Here \(h_{ijk}\) and \(h_{ijkl}\) denote the first and second covariant derivatives of the second
fundamental form \(\vec h\); all sums over indices are ordered. 
The Codazzi and Gauss equations give
\begin{equation}\label{eq:2.0} 
 h_{ijk}=h_{ikj},\qquad
 R_{ijkl}=h_{ik}h_{jl}-h_{il}h_{jk},\qquad
 R_{ij}=Hh_{ij}-\sum_k h_{ik}h_{kj}.
\end{equation}
Defining   \(\Lop=\Delta-\langle X,\nabla\cdot\rangle\), the following equations can be found in \cite{CLW1} if $S$ is constant:
\begin{equation}\label{eq:2.1}
\begin{aligned}
0&=\sum_{i,j,k}h_{ijk}^2+S(1-S),\\
\tfrac13\Lop f_3&=(1-S)f_3+2C,\qquad
\tfrac14\Lop f_4=(1-S)f_4+2A+B,\\
&\sum_{i,j,k,l}h_{ijkl}^2
=S(S-1)(S-2)+3(A-2B).
\end{aligned}
\end{equation}
Henceforth suppose \(S>1\) and put
\[tS=S-1.\]
In order to prove the second gap, we shall suppose $0<t\leq 3/8$. In fact, $t>3/8$ is equivalent to $S>8/5$.
We will use the same notation and symbols as in  \cite{CLW1}.  Since 
\[
 Sf=\sum_{i<j}\lambda_i^2\lambda_j^2(\lambda_i-\lambda_j)^2,
\]
and
\[
 A+2B=\frac13\sum_{ijk}(\lambda_i+\lambda_j+\lambda_k)^2h_{ijk}^2,
\]
we have 
\begin{equation}\label{fidentity}
f\ge0,\qquad Sf=\sum_{i<j}\lambda_i^2\lambda_j^2(\lambda_i-\lambda_j)^2,\qquad A+2B\ge0.
\end{equation}
The following lemma has been proved in Cheng and Wei \cite{CW}:
\begin{lemma}\label{lem:2.2} For a self-shrinker with $S=$constant, 
\begin{equation}\label{eq:3.1}
\begin{aligned}
& \sum_kD_k^2\le \frac{1+2\alpha}3tS^2f,\\
 &C^2 \le \biggl[\frac 13(A+2B)
-\frac {4}3\sum_k\frac {D_k^2}{S+2\lambda_k^2}\biggl]tS^2,
\end{aligned}
\end{equation}
where  $\alpha=\dfrac{\sum_ih_{iii}^2 }{tS^2}\leq 1$. 
\end{lemma}

\section{Key lemmas}

Let $\vec e_1=(\mu_1,\ldots,\mu_n)^T$ and $\vec b=(\mu^2_1-z\mu_1, \cdots, \mu_n^2-z\mu_n)^T$ with $\sqrt S\mu_i=\lambda_i$ and $z=\frac{f_3}{S^{\frac32}}$. Then
\[
 \langle \vec e_1, \vec b\rangle=0,\qquad |\vec b|^2=\frac{ f}{S^2}.
\]
Define
\[
 \vec \psi=(\mu_1^3-cz\mu_1^2, \cdots, \mu_n^3-cz\mu_n^2)^T-(\dfrac{f_4}{S^2}-c\dfrac{f_3^2}{S^3})\vec e_1.
 \]
 If $\vec \psi\neq 0$, define
 \[
  \vec b^* =\vec b-\dfrac{\langle \vec \psi, \vec b\rangle}{|\vec\psi|^2}\vec \psi.
\]
If $\vec \psi =0$, put $\vec b^*=\vec b$.
Set  $f^*=|\vec b^*|^2$, $\rho=S^2\dfrac{f^*}{f}$ when $f>0$, and $\rho=0$ otherwise.
Thus, from $|\vec b^*|^2\leq |\vec b|^2=\dfrac{f}{S^2}$,  $0\le\rho\le1$ and $f^* =\rho\dfrac{f}{S^2}$.
Since $S$ is constant, we have 
$$
\sum\lambda_ih_{iik}=0, \ \text{\rm for any} \ k. 
$$
Putting  $\vec a_k=\frac{1}{\sqrt tS}(h_{11k}, \cdots, h_{nnk})^T$, if 
$$
\sum_i(\lambda_i^3-cz\sqrt S\lambda_i^2)h_{iik}=0,  \ \text{\rm for any} \ k. 
$$
 we have
\begin{equation}\label{eq:3.4}
\langle \vec a_k, \vec e_1\rangle =0, \quad  \langle \vec a_k, \vec \psi \rangle =0, \quad 
u_k=\langle \vec b^*, \vec a_k\rangle=\frac{D_k}{\sqrt tS^2}.
\end{equation}

\begin{lemma}\label{lem:3.1}
\begin{equation}\label{eq:3.2}
2Sf+\left(1+\frac{2(c-1)^2\rho}{1+\rho}\right)f_3^2\le S^3.
\end{equation}
\end{lemma}
\begin{proof}
Suppose first that $f>0$. Putting  $\vec e_2=\frac{S}{\sqrt f}\vec b$, we consider 
the matrix $\Lambda=\operatorname{diag}(\mu_i)$
$$
\Lambda \vec e_1=z\vec e_1+\dfrac{\sqrt f}{S}\vec e_2, \quad \Lambda \vec e_2=\dfrac{\sqrt f}{S}\vec e_1+a\vec e_2+\vec \beta.
$$
Relative to   $\vec e_1$, $\vec e_2$ and their orthogonal complement, the matrix $\Lambda$ can be written 
\[
\Lambda=\begin{pmatrix}
z& \frac{\sqrt f}S&0\\
 \frac{\sqrt f}S&a&\vec \beta^T\\
0&\vec \beta&J_0
\end{pmatrix}.
\]
Since $|\Lambda|^2=1$, we have 
\begin{equation}\label{spectraldefect}
1-2 \dfrac{ f}{S^2}-z^2=a^2+2|\vec\beta|^2+|J_0|^2\ge0.
\end{equation}
Since $(\mu_1^2, \cdots, \mu_n^2)^T=z\vec e_1+\frac{\sqrt{f}}S\vec e_2$, 
$$
(\mu_1^3, \cdots, \mu_n^3)^T=\Lambda (z\vec e_1+\frac{\sqrt{f}}S\vec e_2)=z\Lambda \vec e_1+\frac{\sqrt{f}}S\Lambda \vec e_2.
$$
we obtain
\[
\vec \psi= \dfrac{\sqrt f}S\bigl(\{a+(1-c)z\}\vec e_2+\vec \beta\bigr).
\]
If $\vec \psi\ne0$, 
\begin{equation}\label{eq:3.4-1}
f^*=|\vec b^*|^2=|\vec b|^2-\dfrac{\langle \vec \psi, \vec b\rangle^2}{|\vec\psi|^2}=\dfrac{f}{S^2}\bigl\{1-\frac{(a+(1-c)z)^2}{(a+(1-c)z)^2+|\vec\beta|^2}\bigl\}.
\end{equation}
The definition of $\rho$ and the above formula  yield
\[
\rho=\frac{|\vec\beta|^2}{(a+(1-c)z)^2+|\vec\beta|^2}.
\]
For $0<\rho<1$, it implies
$|\vec\beta|^2=\rho(a+(1-c)z)^2/(1-\rho)$. Hence, we infer
\begin{align*}
a^2+2|\vec\beta|^2=a^2+\frac{2\rho}{1-\rho}(a+(1-c)z)^2
={}&\frac{1+\rho}{1-\rho}
 \left(a+\frac{2\rho}{1+\rho}(1-c)z\right)^2
 +\frac{2\rho}{1+\rho}(1-c)^2z^2.
\end{align*}
In view of  \eqref{spectraldefect}, we have 
$$
1\geq 2 \dfrac{ f}{S^2}+\bigl\{1+\frac{2\rho}{1+\rho}(1-c)^2\bigl\}z^2.
$$
The claimed estimate for $0<\rho<1$ is obtained.
If $\rho=1$, including the case $\vec \psi=0$, then $a=(c-1)z$.
\eqref{spectraldefect} yields the same assertion.
For $\rho=0$, the assertion follows from \eqref{spectraldefect}.
If $f=0$, our convention is $\rho=0$, and
\eqref{eq:3.2} reduces to $f_3^2\le S^3$, which is obvious.
\end{proof}
Define
\[
\mathcal G_{kl}=\frac{1}{tS^2}\sum_{i,j}h_{ijk}h_{ijl},\qquad
\mathcal M_{kl}=4\delta_{kl}+12\mathcal G_{kl}.
\]
Then $\mathcal G=(\mathcal G_{kl})$ is a positive semi-definite matrix and $\tr\mathcal G=1$, and $\mathcal M=(\mathcal M_{kl})$ is
positive definite.
\begin{lemma}\label{lem:3.2}
For any vector $\vec w=(w_1, \cdots, w_n)^T$, we have 
\begin{equation}\label{eq:3.5}
\sum_{k,l}\mathcal G_{kl}w_kw_l+2\sum_{k,l}\mathcal G_{kl}w_k\mu_kw_l\mu_l\leq |\vec w|^2,
\end{equation}
Furthermore, if $\sum_i(\lambda_i^3-cz\sqrt S\lambda_i^2)h_{iik}=0,  \ \text{\rm for any} \ k$,
\begin{equation}\label{eq:3.6}
\sum_kD_k^2+2\sum_k\mu_k^2D_k^2\le tS^4f^*,\qquad
(\sum_ku_kw_k)^2\leq f^*\sum_{k,l}\mathcal G_{kl}w_kw_l.
\end{equation}
If $f^*>0$ and $\vec v = (\mu_1u_1, \cdots, \mu_nu_n)^T \ne0$, then
\begin{equation}\label{eq:3.7}
\vec v^T\mathcal M^{-1}\vec v
\ge\frac{f^*\sum_{k}\lambda_k^2D_k^2}{16tS^5f^*-24\sum_{k}\lambda_k^2D_k^2}.
\end{equation}
If $f^*=0$, then $\vec v=0$.
\end{lemma}
\begin{proof}
Let $\vec e_a$ be an orthonormal eigenbasis of $\mathcal G$, with
eigenvalues $g_a\ge0$. For $\vec w\in\R^n$, write $T_w=(\frac1{\sqrt tS}\sum_kh_{ijk}w_k)$ and 
$T(w, e_a, e_b)=\frac1{\sqrt tS}\sum_{i,j,k}h_{ijk}w_k(e_a)_i(e_b)_j$, where $(e_a)_i$ denotes
the $i$th component  of the vector $\vec e_a$.
The symmetric matrices
\[
\Lambda,\qquad T_{e_b}/\sqrt{g_b}\quad(g_b>0)
\]
are orthonormal since $S$ is constant.
Indeed, $\langle\Lambda,T_w\rangle=\frac1{\sqrt tS}\sum_{i,k}\mu_ih_{iik}w_k=0$ for every $\vec w$, while
$\langle T_{e_a},T_{e_b}\rangle=g_a\delta_{ab}$ and
$|\Lambda|^2=1$. Set
$B_a(w)=(\vec w\otimes \vec e_a+\vec e_a\otimes \vec w)/2$. Then
\[
|B_a(w)|^2=\tfrac12|\vec w|^2+\tfrac12\langle \vec w, \vec e_a\rangle^2,
\quad \langle B_a(w),\Lambda\rangle=\langle\Lambda w,e_a\rangle,
\quad \langle B_a(w),T_{e_b}\rangle=T(w,e_a,e_b).
\]
Putting
$$
R_a(w)=B_a(w)-\langle B_a(w),\Lambda\rangle\Lambda-\sum_{b, g_b>0}\langle B_a(w), T_{e_b}\rangle \dfrac{T_{e_b}}{g_b},
$$
we have 
$$
|R_a(w)|^2=|B_a(w)|^2-\langle B_a(w),\Lambda\rangle^2-\sum_{b, g_b>0}\langle B_a(w), T_{e_b}\rangle^2 \dfrac{1}{g_b}.
$$
By $g_a$ multiplying  the above formula  and taking summation, we get   
 $$
0\leq \sum_ag_a|R_a(w)|^2=\sum_ag_a|B_a(w)|^2-\sum_ag_a\langle B_a(w),
\Lambda\rangle^2-\sum_{a, b, g_b>0}\dfrac{g_a}{g_b}\langle B_a(w), T_{e_b}\rangle^2.
$$
 According to 
 $$
 \sum_ag_a\langle B_a(w),\Lambda\rangle^2=\langle \mathcal G(\Lambda \vec w),\Lambda\vec w\rangle
 =\vec w^T\Lambda\mathcal G\Lambda\vec w
 $$
 and 
 $$
 \sum_ag_a|B_a(w)|^2=\dfrac12|\vec w|^2\sum_ag_a+\dfrac12\sum_ag_a\langle \vec w, \vec e_a\rangle^2
 =\dfrac12|\vec w|^2+\dfrac12 \vec w^T\mathcal G\vec w,
 $$
we conclude
\begin{equation}\label{expandedbessel}
\vec w^T\Lambda\mathcal G\Lambda \vec w+
\sum_{a,b:g_a,g_b>0}\frac{g_a}{g_b}T(w,e_a,e_b)^2
\le\tfrac12|\vec w|^2+\tfrac12\vec w^T\mathcal G \vec w.
\end{equation}
If $g_a=0$, then $T_{e_a}=0$; all terms containing that index
vanish. The symmetry in $a,b$ and   
$g_a/g_b+g_b/g_a\ge2$  give
\[
\sum_{a,b:g_a,g_b>0}\frac{g_a}{g_b}T(w,e_a,e_b)^2
\ge\sum_{a,b}T(w,e_a,e_b)^2=|T_w|^2=\vec w^T\mathcal G \vec w,
\]
where we used that  $\{\vec e_a\}$ is an orthonormal basis.
Substitution in \eqref{expandedbessel} gives
\[
\vec w^T\mathcal G \vec w+\vec w^T\Lambda\mathcal G\Lambda \vec w
\le\frac12|\vec w|^2+\frac12\vec w^T\mathcal G \vec w.
\]
This proves (\ref{eq:3.5}).
\\
Let $E_i$ denote the $i$th
standard coordinate vector. Define $s_i=\dfrac{h_{iii}}{\sqrt t S}$, since $\langle \vec a_{k}, \vec e_1\rangle=0$ from (\ref{eq:3.4}),
\[
 s_i=\langle \vec a_i, E_i-\mu_i\vec e_1\rangle,
 \qquad |E_i-\mu_i\vec e_1|^2=1-\mu_i^2.
\]
Cauchy-Schwarz inequality therefore gives
$s_i^2\le(1-\mu_i^2)|\vec a_i|^2$.
Since $h_{ijk}$ are symmetric for $i,\ j, \ k $, $\sum_{i,j,k}h_{ijk}^2=tS^2$  gives
\[
1=\sum_i\bigl(3|\vec a_i|^2-2s_i^2\bigr)
 +6\dfrac{1}{tS^2}\sum_{i<j<k}h_{ijk}^2
\ge\sum_i(1+2\mu_i^2)|\vec a_i|^2.
\]
Together with (\ref{eq:3.4}), this implies
$$
\sum_i(1+2\mu_i^2)u_i^2\le \sum_i(1+2\mu_i^2)|\vec a_i|^2f^*\le f^*.
$$
According to (\ref{eq:3.4}) and $\sum_i(\lambda_i^3-cz\sqrt S\lambda_i^2)h_{iik}=0,  \ \text{\rm for any} \ k$, we have 
\begin{equation}\label{eq:3.9}
u_k=\langle \vec b^*, \vec a_k\rangle=\langle \vec b, \vec a_k\rangle=\dfrac{1}{\sqrt tS^2}D_k
\end{equation}
$$
\sum_k(1+2\mu_k^2)D_k^2\le tS^4 f^*.
$$
Also, for every $\vec w$,
\[
\sum_ku_kw_k=\sum_ib_i^*\sum_k\dfrac{h_{iik}}{\sqrt t S}w_k,
\qquad
(\sum_ku_kw_k)^2\le f^*\sum_{k,l}\mathcal G_{kl}w_k w_l.
\]
This proves \eqref{eq:3.6}. 
\\
For $f^*>0$, for any vector $\vec w$, from (\ref{eq:3.5}) we infer 
\[
\sum_{k,l}\mathcal M_{kl}w_kw_l=4|\vec w|^2+12\sum_{k,l}\mathcal G_{kl}w_kw_l
\leq 16|\vec w|^2-24\sum_{k,l}\mathcal G_{kl}w_k\mu_kw_l\mu_l.
\]
The second formula in \eqref{eq:3.6} gives
\[
\sum_{k,l}\mathcal M_{kl}w_kw_l
\leq 16|\vec w|^2-\dfrac{24}{f^*}(\sum_{k}u_k\mu_kw_k)^2.
\]
Putting  $\vec v=(\mu_1u_1, \cdots, \mu_nu_n)^T$,  (\ref{eq:3.9}) yields 
$$
|\vec v|^2=\sum_k\mu_k^2u_k^2= \dfrac{1}{tS^5}\sum_{k}\lambda_k^2D_k^2.
$$
Taking $\vec w=\vec v$,  the above inequalities   give
\[
\begin{aligned}
0&<v^T\mathcal Mv\le16\dfrac{1}{tS^5}\sum_{k}\lambda_k^2D_k^2
-\frac{24\biggl(\dfrac{1}{tS^5}\sum_{k}\lambda_k^2D_k^2\biggl)^2}{f^*}\\
&=\dfrac{1}{tS^5}\sum_{k}\lambda_k^2D_k^2\left(16-\frac{\dfrac{24}{tS^5}\sum_{k}\lambda_k^2D_k^2}{f^*}\right).
\end{aligned}
\]
when $\vec v\ne0$. Since $\mathcal M$ is positive definite,  we have 
$$
|\vec v|^4=\langle \vec v, \vec v\rangle^2=\langle \mathcal M^{1/2} \vec v, \mathcal M^{-1/2} \vec v\rangle^2
\leq \vec v^T\mathcal M^{-1}\vec v\cdot \vec v^T\mathcal M\vec v.
$$
Hence, we obtain
\[
\vec v^T\mathcal M^{-1}\vec v\ge\frac{|\vec v|^4}{\vec v^T\mathcal M \vec v}
\ge\frac{\dfrac{f^*}{tS^5}\sum_{k}\lambda_k^2D_k^2}{16f^*-24\dfrac{1}{tS^5}\sum_{k}\lambda_k^2D_k^2}.
\]
For $\vec v=0$, the corresponding lower bound is zero. 
\end{proof}
Recall that $b_i=\mu_i^2-z\mu_i$ are the components of $\vec b$.
We have
\[
 \sum_i b_i^2=\frac{f}{S^2},\qquad
 \sum_i\mu_i^4=\frac{f_4}{S^2}=\frac{f}{S^2}+z^2,
 \qquad 0\le z^2\le1.
\]
By Lemma~\ref{lem:3.1}, $2f/S^2+z^2\le1$.

\begin{lemma}\label{commutator}
Let $z\ge0$ and $c\ge1$. Suppose that, for each $k$,
\[
 \sum_i\lambda_i h_{iik}=0,\qquad
 \sum_i(\lambda_i^3-cz\sqrt S\lambda_i^2)h_{iik}=0.
\]
Then
\[
4(A-B)\leq 2tS^3.
\]
\end{lemma}
\begin{proof}
Put $T_{ijk}=h_{ijk}/(\sqrt t S)$ and $\beta=cz$. Then
\[
 \sum_{ijk}T_{ijk}^2=1,\qquad
 \sum_i\mu_iT_{iik}=0,\qquad
 \sum_i\mu_i^2(\mu_i-\beta)T_{iik}=0.
\]
Since $T_{ijk}$ is symmetric in its three indices,
\[
 \frac{4(A-B)}{tS^3}
 =2\sum_{ijk}(\mu_i-\mu_j)^2T_{ijk}^2.
\]
For three distinct indices $i,j,k$,
\begin{align*}
 &(\mu_i-\mu_j)^2+(\mu_j-\mu_k)^2+(\mu_k-\mu_i)^2\\
 &\qquad=3(\mu_i^2+\mu_j^2+\mu_k^2)
             -(\mu_i+\mu_j+\mu_k)^2\le3.
\end{align*}
By the symmetry of $T_{ijk}$, summing over the six permutations gives
\[
2\sum_{\text{permutations of }(i,j,k)}(\mu_i-\mu_j)^2T_{ijk}^2
\le12T_{ijk}^2.
\]
The three permutations of $(i,i,k)$ contribute
$4(\mu_i-\mu_k)^2T_{iik}^2$, whereas $(k,k,k)$ contributes zero. It remains to prove, for each $k$,
\begin{equation}\label{repeated}
\sum_{i\ne k} [2(\mu_i-\mu_k)^2-3]T_{iik}^2\leq T_{kkk}^2.
\end{equation}
If every coefficient on the left is nonpositive, \eqref{repeated} is immediate.
Otherwise, there is at most one positive coefficient. Indeed,
\[
 (\mu_i-\mu_k)^2+(\mu_j-\mu_k)^2
 \le3(\mu_i^2+\mu_j^2+\mu_k^2)\le3.
\]
The first inequality follows by expanding: the difference between its
right and left sides is
$(\mu_i-\mu_j)^2+(\mu_i+\mu_j+\mu_k)^2$. If such an index $i$ exists, put $P=2(\mu_i-\mu_k)^2-3>0$. All other coefficients are at most $-P$. Choose $(Q_j)$ such that $\sum_j Q_jT_{jjk}=0$ and
\begin{equation}\label{Qcondition}
P Q_k^2+\sum_{j\ne i,k}Q_j^2\leq Q_i^2,\qquad Q_i\ne0.
\end{equation}
Then
\[
 Q_iT_{iik}=-Q_kT_{kkk}-\sum_{j\ne i,k}Q_jT_{jjk}.
\]
By Cauchy's inequality,
\[
 Q_i^2T_{iik}^2\le
 \left(Q_k^2+\frac1P\sum_{j\ne i,k}Q_j^2\right)
 \left(T_{kkk}^2+P\sum_{j\ne i,k}T_{jjk}^2\right).
\]
Thus \eqref{Qcondition} implies
$PT_{iik}^2\le T_{kkk}^2+P\sum_{j\ne i,k}T_{jjk}^2$,
and hence \eqref{repeated} holds.
Put $v=\mu_k$, $u=\mu_i$, and $L=\sum_{j\ne i,k}\mu_j^2$.
Since $u^2+v^2+L=1$ and $P=2(u-v)^2-3$, we have
\[
u^2-L=P(1-v^2)+2(1-v^2+uv)^2.
\]
If $v^2\leq1/2$, then
\[
 u^2-L\ge P(1-v^2)\ge Pv^2.
\]
Thus $Q_j=\mu_j$ satisfies \eqref{Qcondition}.
Now suppose $v^2>1/2$. If $v<0$, then
\[
 z=\sum_j\mu_j^3\le -|v|^3+(1-v^2)^{3/2}<0,
\]
a contradiction. Thus $v>0$. Also $P>0$ implies $u<0$, since $u\ge0$ would give $(u-v)^2\le u^2+v^2\le1$. Moreover, $u^2\ge L$, and hence $|\mu_j|\le -u$ for $j\ne k$.
Consequently $\mu_j^3\ge u\mu_j^2$ for every $j\ne k$, and hence
\[
P(1-v^2)^2+L^2\leq u^4,
\qquad \sqrt L\leq -u,
\qquad \beta\geq z\geq v^3+u(1-v^2).
\]
For the first inequality, multiply $u^2-L\ge P(1-v^2)$ by
$u^2+L=1-v^2$ to obtain $u^4-L^2\ge P(1-v^2)^2$. In this case take
\[
Q_j=\mu_j^2(\mu_j-\beta)+v(\beta-v)_+\mu_j.
\]
The two gradient relations give $\sum_jQ_jT_{jjk}=0$.
If $\beta\ge v$, then $Q_k=0$. If $\beta<v$, then
\[
 v^2(v-\beta)\le(1-v^2)(\beta-u)
 \quad\Longleftrightarrow\quad
 v^3+u(1-v^2)\le\beta,
\]
which was proved above. Thus
$0\leq Q_k=v^2(v-\beta)_+\leq(1-v^2)(\beta-u)$, and $Q_i<0$. Since $|\mu_j|\leq -u$ for $j\ne k$, the triangle inequality gives
\begin{align*}
(P Q_k^2+\sum_{j\ne i,k}Q_j^2)^{1/2}
&\leq(\beta-u)\sqrt{P(1-v^2)^2+L^2}+v(\beta-v)_+\sqrt L\\
&\leq(\beta-u)u^2-uv(\beta-v)_+=|Q_i|.
\end{align*}
Thus \eqref{Qcondition} holds. Summing \eqref{repeated} together with the distinct-index estimate proves the assertion.
\end{proof}

\begin{lemma}\label{spectral2}
Suppose $f>0$ and let $r=\frac{S^2}{f}\sum_i\mu_i^3b_i$. We have
\begin{equation}\label{jbound}
\frac{f}{2S^2}-\sum_i\mu_i^2b_i^2\geq\frac{S^2}{S^2-f_4}\left(\sum_i\mu_i^3b_i\right)^2.
\end{equation}
Moreover, for $a_0,b_0\ge0$,
\begin{equation}\label{projectbound}
\begin{aligned}
 &a_0^2f^*-b_0^2\left(\frac{f}{2S^2}-\sum_i\mu_i^2b_i^2\right)\\
 &\qquad\le\frac{f}{S^2}(a_0-b_0|r-cz|)_+^2
 -\frac{b_0^2f}{S^2}\left(\frac{f_4}{S^2}-\frac12\right).
\end{aligned}
\end{equation}
\end{lemma}
\begin{proof}
By Lemma~\ref{lem:3.1}, $S^2-f_4\ge f>0$. For $s\in\mathbb R$, the off-diagonal squared norm of the matrix with entries $\mu_i b_j+\mu_j b_i-s\mu_i\mu_j$ is
\[
\frac{2f}{S^2}-4\sum_i\mu_i^2b_i^2+4s\left(\sum_i\mu_i^3b_i\right)+s^2(1-\frac{f_4}{S^2})\geq0.
\]
Indeed, the full squared norm is $\frac{2f}{S^2}+s^2$, since
$\langle\vec e_1,\vec b\rangle=0$ and $|\vec e_1|=1$.
The sum of its diagonal squares is
$4\sum_i\mu_i^2b_i^2-4s\left(\sum_i\mu_i^3b_i\right)+s^2(\frac{f_4}{S^2})$. Subtracting gives the displayed expression.
Taking $s=-\dfrac{2S^2}{S^2-f_4}\sum_i\mu_i^3b_i$ gives
\[
 0\le\frac{2f}{S^2}-4\sum_i\mu_i^2b_i^2-\frac{4S^2}{S^2-f_4}\left(\sum_i\mu_i^3b_i\right)^2,
\]
which is \eqref{jbound}.
In the notation of the proof of Lemma~\ref{lem:3.1}, $r=z+a$ and
\begin{align*}
 \sum_i\mu_i^2b_i^2
 &=\frac{f}{S^2}\left(\frac{f}{S^2}+a^2+|\vec\beta|^2\right),\\
 \frac{f}{2S^2}-\sum_i\mu_i^2b_i^2
 &=\frac{f}{S^2}\left(\frac{f_4}{S^2}-\frac12
                            +|\vec\beta|^2+|J_0|^2\right).
\end{align*}
Let $P_0$ denote the orthogonal projection onto $\{\vec e_1,\vec\psi\}^\perp$. Since
$\vec\psi=(r-cz)\vec b+\frac{\sqrt f}{S}\vec\beta$, its projection is zero.
For $|\tau|\le b_0$, the identities $P_0\vec b=\vec b^*$ and $P_0\vec\psi=0$ give
\[
 a_0\vec b^*=P_0\bigl((a_0-\tau(r-cz))\vec b
                      -\tau\frac{\sqrt f}{S}\vec\beta\bigr).
\]
Since $P_0$ is an orthogonal projection and $\vec b\perp\vec\beta$,
\[
a_0^2f^*\leq \frac{f}{S^2}(a_0-\tau(r-cz))^2+\frac{b_0^2f}{S^2}|\vec\beta|^2.
\]
Minimizing over $|\tau|\le b_0$, we obtain
\[
a_0^2f^*\le\frac{f}{S^2}(a_0-b_0|r-cz|)_+^2
 +\frac{b_0^2f}{S^2}|\vec\beta|^2.
\]
Subtracting $b_0^2(\frac{f}{2S^2}-\sum_i\mu_i^2b_i^2)$ and using its expression above yields
\begin{align*}
 a_0^2f^*-b_0^2(\frac{f}{2S^2}-\sum_i\mu_i^2b_i^2)
 &\le \frac{f}{S^2}(a_0-b_0|r-cz|)_+^2
       -\frac{b_0^2f}{S^2}(\frac{f_4}{S^2}-\frac12+|J_0|^2)\\
 &\le \frac{f}{S^2}(a_0-b_0|r-cz|)_+^2-\frac{b_0^2f}{S^2}(\frac{f_4}{S^2}-\frac12).
\end{align*}
This proves \eqref{projectbound}.
\end{proof}
\section{Main estimates}\label{tensorsection}

Defining 
\[
a_{ij}=\sum_kh_{ik}h_{kj}-\frac{f_3}{S}h_{ij},\qquad \sum_{i,j}a_{ij}^2=f.
\]
$u_{ijkl}$ introduced by Cheng and Wei \cite{CW} is defined by
\[
u_{ijkl}=\tfrac14(h_{ijkl}+h_{jkli}+h_{klij}+h_{lijk}).
\]
The following formulas can be found in Cheng, Li and Wei \cite{CLW1}:
\begin{equation*}
\sum_{ijkl}h_{ijkl}^2=\sum_{i,j,k,l}u_{ijkl}^2+\frac32Sf.
\end{equation*}
and
\begin{equation}\label{eq:4.3}
\sum_{ijkl}u_{ijkl}^2=t(2t-1)S^3+3(A-2B)-\frac32Sf.
\end{equation}
Differentiating $S$ twice and using the Ricci commutation identities, as in \cite{CLW1}, gives
\begin{equation}\label{newcontractions}
\sum_{i,j}u_{ijkl}h_{ij}=-\sum_{i,j}h_{ijk}h_{ijl}+\frac S2\sum_{i}h_{ki}h_{il}-\frac{f_3}2h_{kl},
\qquad
\sum_{i,j,k}u_{ijkl}h_{ijk}=-\frac32\lambda_lD_l.
\end{equation}
\[
\sum_i\lambda_i h_{iikl}=-\sum_{i,j}h_{ijk}h_{ijl}.
\]
Define, for $\vec \xi=(\xi_1, \cdots, \xi_n)^T$,
\begin{equation}
\begin{aligned}
 &V_{ijkl}(\xi)=\xi_ih_{jkl}+\xi_jh_{kli}+\xi_kh_{lij}+\xi_lh_{ijk}, \\
 & X_{ijkl}=h_{ij}a_{kl}+a_{ij}h_{kl},
  \quad Z_{ijkl}=h_{ij}h_{kl}.
\end{aligned}
 \end{equation}
Since $S$ is constant, we obtain the following by a direct computation:
\begin{lemma}\label{lem:4.1} For a self-shrinker with constant $S$, we have 
$$
\begin{aligned}
& \sum_{i,j,k,l}X_{ijkl}^2=2S(f_4-\dfrac{f_3^2}S)=2Sf, \quad 
 \sum_{i,j,k,l}Z_{ijkl}^2=S^2,\\
&\sum_{i,j,k,l}X_{ijkl}V_{ijkl}(\xi)=4\sum_l\xi_l\lambda_lD_l, \quad \sum_{i,j,k,l}V_{ijkl}(\xi)Z_{ijkl}=0,\\
& \sum_{i,j,k,l}X_{ijkl}u_{ijkl}=-2A+Sf+\dfrac{2f_3C}S,  \quad  \sum_{i,j,k,l}u_{ijkl}Z_{ijkl}=-C,\\
 &\sum_{i,j,k,l}u_{ijkl}V_{ijkl}(\xi)=-6\sum_l\xi_l\lambda_lD_l,\quad  \sum_{i,j,k,l}X_{ijkl}Z_{ijkl}=0,\\
& \sum_{i,j,k,l}V_{ijkl}^2(\xi)=tS^2\bigl(4|\vec \xi|^2+12\sum_{i,j}\mathcal G_{ij}\xi_i\xi_j\bigl).
\end{aligned}
$$
\end{lemma}
\begin{proof}
In the principal frame,
$a_{ij}=(\lambda_i^2-f_3\lambda_i/S)\delta_{ij}$, and hence
\[
 \sum_{ij}h_{ij}a_{ij}=f_3-\frac{f_3}{S}S=0,
 \qquad \sum_{ij}a_{ij}^2=f.
\]
Expanding the products in $X=h\otimes a+a\otimes h$ and $Z=h\otimes h$ gives
\[
 |X|^2=2Sf+2\left(\sum_{ij}h_{ij}a_{ij}\right)^2=2Sf,
 \qquad |Z|^2=S^2,\qquad \langle X,Z\rangle=0.
\]
Since $\sum_i\lambda_i h_{iik}=0$, the terms in $\langle X,V(\xi)\rangle$
containing this contraction vanish. The other four terms each equal
$\sum_k\xi_k\lambda_kD_k$. Thus
\[
 \langle X,V(\xi)\rangle=4\sum_k\xi_k\lambda_kD_k,
 \qquad\langle Z,V(\xi)\rangle=0.
\]
Using the first identity in \eqref{newcontractions}, we obtain
\begin{align*}
 \langle u,X\rangle
 &=2\sum_{kl}a_{kl}\left(-\sum_{ij}h_{ijk}h_{ijl}
                              +\frac S2a_{kl}\right)\\
 &=-2\sum_{ijk}\left(\lambda_k^2-\frac{f_3}{S}\lambda_k\right)h_{ijk}^2+Sf
 =-2A+\frac{2f_3C}{S}+Sf,\\
 \langle u,Z\rangle
 &=\sum_{kl}h_{kl}\left(-\sum_{ij}h_{ijk}h_{ijl}
                              +\frac S2a_{kl}\right)=-C.
\end{align*}
The full symmetry of $u$ and the second identity in
\eqref{newcontractions} give
\[
 \langle u,V(\xi)\rangle
 =4\sum_l\xi_l\sum_{ijk}u_{ijkl}h_{ijk}
 =-6\sum_l\xi_l\lambda_lD_l.
\]
By the symmetry of $h_{ijk}$, the six cross terms in the expansion of
$|V(\xi)|^2$ are equal. Hence
\[
 |V(\xi)|^2=4tS^2|\vec\xi|^2
       +12\sum_{i,j,k,l}h_{ijk}h_{ijl}\xi_k\xi_l
 =tS^2\left(4|\vec\xi|^2
       +12\sum_{k,l}\mathcal G_{kl}\xi_k\xi_l\right).
\]
\end{proof}

Let $\bar X$ and $\bar Z$ be the symmetrizations of $X$ and $Z$, respectively:
\begin{align*}
3\bar X_{ijkl}&=h_{ij}a_{kl}+h_{ik}a_{jl}+h_{il}a_{jk}
+a_{ij}h_{kl}+a_{ik}h_{jl}+a_{il}h_{jk},\\
3\bar Z_{ijkl}&=h_{ij}h_{kl}+h_{ik}h_{jl}+h_{il}h_{jk}.
\end{align*}
By the symmetry of $u$ and $V(\xi)$,
\[
 \langle u,\bar X\rangle=\langle u,X\rangle,\qquad
 \langle V(\xi),\bar X\rangle=\langle V(\xi),X\rangle,
\]
and similarly for $Z$ and $\bar Z$. In a frame diagonalizing $h$, the nonzero components are, up to permutation,
\[
 \bar Z_{iiii}=\lambda_i^2,\qquad
 \bar Z_{iijj}=\frac13\lambda_i\lambda_j,
 \qquad \bar X_{iiii}=2\lambda_i a_{ii},\qquad
 \bar X_{iijj}=\frac13(\lambda_i a_{jj}+\lambda_j a_{ii}).
\]
Summing over the indices and using $\sum_i\lambda_i a_{ii}=0$, we obtain
\begin{align*}
 |\bar Z|^2
 &=\sum_i\lambda_i^4+\frac23\sum_{i<j}\lambda_i^2\lambda_j^2
   =\frac13(S^2+2f_4),\\
 |\bar X|^2
 &=4\sum_i\lambda_i^2a_{ii}^2
   +\frac23\sum_{i<j}(\lambda_i a_{jj}+\lambda_j a_{ii})^2\\
 &=\frac23Sf+\frac83\sum_i\lambda_i^2a_{ii}^2,\\
 \langle\bar X,\bar Z\rangle
 &=2\sum_i\lambda_i^3a_{ii}
   +\frac23\sum_{i<j}\lambda_i\lambda_j
                        (\lambda_i a_{jj}+\lambda_j a_{ii})
   =\frac43\sum_i\lambda_i^3a_{ii}.
\end{align*}
Since $a_{ii}=Sb_i$ and $\lambda_i=\sqrt S\mu_i$, these identities
can be written as
\begin{equation}\label{symnorms}
\begin{aligned}
 |\bar Z|^2&=\frac{S^2+2f_4}{3},\\
 |\bar X|^2&=\frac23Sf+\frac83S^3\sum_i\mu_i^2b_i^2,\\
 \langle\bar X,\bar Z\rangle
 &=\frac43\sum_i\lambda_i^3a_{ii}
 =\frac43S^{5/2}\sum_i\mu_i^3b_i.
\end{aligned}
\end{equation}

Since $S$ is constant, $|\lambda_i|\le\sqrt S$ and $|H|\le\sqrt{nS}$.
The first identity in \eqref{eq:2.1} gives
$\sum h_{ijk}^2=S(S-1)$. Moreover,
\[
 |A|\le S\sum h_{ijk}^2,\qquad
 |B|\le S\sum h_{ijk}^2.
\]
It follows from \eqref{eq:2.1} that $\sum h_{ijkl}^2$ is also bounded.
Moreover, the Gauss equation gives
\[
 \operatorname{Ric}(e_i,e_i)=H\lambda_i-\lambda_i^2
 \ge-(\sqrt n+1)S,
\]
so the Ricci curvature is bounded below. Let
\[
F=\frac S4f_4-\frac{137}{600}f_3^2+\frac{f_3^4}{120S^3}.
\]
Since $F$ is bounded, the generalized maximum principle of Cheng and Peng \cite{CP} gives a sequence $\{p_m\}$ such that
\[
 F(p_m)\to\sup_M F,\qquad |\nabla F|(p_m)\to0,
 \qquad\limsup_m\mathcal L F(p_m)\le0.
\]
At each $p_m$, choose a principal frame and a unit normal with $f_3(p_m)\ge0$.
Taking a subsequence, we may assume that $h_{ij},h_{ijk},h_{ijkl}$ converge. We use the same notation for their limits.
Put
\[
c=\frac{137}{100}-\frac{z^2}{10},\qquad d=\frac{137}{100}-\frac{3z^2}{10}.
\]
where $z=f_3/S^{3/2}$ at the limit. Differentiating $F$, we have
\begin{align*}
 \nabla_k F
 &=S\sum_i\lambda_i^3h_{iik}
   -\left(\frac{137}{300}f_3-\frac{f_3^3}{30S^3}\right)\nabla_k f_3\\
 &=S\sum_i\lambda_i^3h_{iik}-c f_3\sum_i\lambda_i^2h_{iik}.
\end{align*}
Together with $\nabla_k S=2\sum_i\lambda_i h_{iik}=0$, this gives the gradient equations
\begin{equation}\label{generalgrad}
\sum_i\lambda_i h_{iik}=0,\qquad
\sum_i(\lambda_i^3-c f_3\lambda_i^2/S)h_{iik}=0.
\end{equation}
Furthermore, $c\ge127/100>1$. The product rule for $\mathcal L$ gives
\[
 \mathcal L(f_3^2)=2f_3\mathcal Lf_3+2|\nabla f_3|^2,
 \quad
 \mathcal L(f_3^4)=4f_3^3\mathcal Lf_3+12f_3^2|\nabla f_3|^2.
\]
Therefore
\begin{align*}
 \mathcal LF
 &=\frac S4\mathcal Lf_4
   -\frac{cf_3}{3}\mathcal Lf_3-\frac d3|\nabla f_3|^2\\
 &=S(1-S)f_4+S(2A+B)
   -cf_3\bigl((1-S)f_3+2C\bigr)-3d\sum_kD_k^2.
\end{align*}
Dividing by $S$ and using $1-S=-tS$, we obtain
\begin{equation}\label{lapF}
0\geq\frac1S\limsup\Lop F
=-tSf_4+2A+B+ctf_3^2-\frac{2cf_3 C}{S}-\frac{3d}{S}\sum_k D_k^2.
\end{equation}

Put $m=C/(tS^{5/2})$. If $f^*>0$, define
\[
x=\frac{\sum_kD_k^2}{tS^4 f^*},\qquad y=\frac{\sum_k\lambda_k^2D_k^2}{tS^5 f^*};
\]
otherwise set $x=y=0$. From Lemma~\ref{lem:3.2}, $x,y\geq0$ and $x+2y\leq1$. In what follows, $x^2/(x+2y)$ is defined as zero at $(0,0)$. Cauchy--Schwarz and \eqref{eq:3.1} imply
\begin{equation}\label{ABlower}
A+2B\geq 3tS^3m^2+4tS^3f^*\frac{x^2}{x+2y}.
\end{equation}
Indeed, $\sum D_k^2/(S+2\lambda_k^2)\geq (\sum D_k^2)^2/(S\sum D_k^2+2\sum\lambda_k^2D_k^2)$, with zero interpreted by continuity.

\begin{proposition}\label{newtensor}
For every $\delta>0$, put $k_\delta=4+8\delta/3$. Under the sense of the above limit, the following inequalities hold:
\begin{align}
0\le\Phi_\delta(t)={}&2t-1+
 \left(k_\delta-\frac{3/2+2\delta-2\delta^2}{t}\right)\frac{f}{S^2}\notag\\
 &+k_\delta(1-c)z^2+(2ck_\delta-4\delta)zm-(15+4\delta)m^2\notag\\
 &-\frac{S^2}{3t(S^2+2f_4)}
                  \left(3tm+4\delta\sum_i\mu_i^3b_i\right)^2\notag\\
 &-\frac{8\delta^2}{3t}
                  \left(\frac{f}{2S^2}-\sum_i\mu_i^2b_i^2\right)\notag\\
 &+f^*\left(3dk_\delta x-\frac{4(5+4\delta/3)x^2}{x+2y}
                  -\frac{(3+2\delta)^2y}{t(4-6y)}\right),\label{Phi}\\
0\le\Psi(t)={}&\frac23t(2t+1)-\frac{f}{S^2}-2tm^2\notag\\
 &-f^*\left(\frac{8tx^2}{3(x+2y)}+\frac{6y}{4-6y}\right).\label{Psi}
\end{align}
\end{proposition}
\begin{proof}
By Lemma~\ref{lem:4.1}, $\bar Z$ is orthogonal to all $V(\xi)$.
Put $\vec w=(\lambda_1D_1,\ldots,\lambda_nD_n)^T$. Expansion of the nonnegative square gives
\begin{align*}
 0\le{}&|u-\delta\bar X+\eta\bar Z+V(\xi)|^2\\
 ={}&|u-\delta\bar X|^2+\eta^2\frac{S^2+2f_4}{3}
       -2\eta\left(C+\frac{4\delta}{3}S^{5/2}\left(\sum_i\mu_i^3b_i\right)\right)\\
 &+tS^2\vec\xi^T\mathcal M\vec\xi
       -(12+8\delta)\langle\vec\xi,\vec w\rangle.
\end{align*}
The minimizing choices are
\[
 \eta=\frac{3C+4\delta S^{5/2}\sum_i\mu_i^3b_i}{S^2+2f_4},
 \qquad \vec\xi=\frac{6+4\delta}{tS^2}\mathcal M^{-1}\vec w.
\]
Completing the squares in $\eta$ and $\vec\xi$ therefore gives
\begin{equation}\label{symprojection}
\begin{aligned}
 |u-\delta\bar X|^2\ge{}&
 \frac{\left(3C+4\delta S^{5/2}\sum_i\mu_i^3b_i\right)^2}
      {3(S^2+2f_4)}\\
 &+\frac{(6+4\delta)^2}{tS^2}\vec w^T\mathcal M^{-1}\vec w.
\end{aligned}
\end{equation}
For $f^*>0$, Lemma~\ref{lem:3.2} and
$\sum_k\lambda_k^2D_k^2=tS^5f^*y$ imply
\[
 \vec w^T\mathcal M^{-1}\vec w
 \ge\frac{tS^5f^*y}{16-24y}.
\]
Consequently the second term in \eqref{symprojection} is at least
$(3+2\delta)^2S^3f^*y/(4-6y)$. Its denominator is positive because
$x+2y\le1$ implies $4-6y\ge1$. When $f^*=0$, both sides of this estimate
are zero. On the other hand, \eqref{eq:4.3}, Lemma~\ref{lem:4.1} and \eqref{symnorms} give
\begin{align*}
|u-\delta\bar X|^2={}&t(2t-1)S^3+3(A-2B)+4\delta A\\
&+(-3/2-2\delta+2\delta^2)Sf
-\frac{8\delta^2}{3}S^3(\frac{f}{2S^2}-\sum_i\mu_i^2b_i^2)-4\delta f_3C/S.
\end{align*}
We use $3(A-2B)+4\delta A=k_\delta(2A+B)-(5+4\delta/3)(A+2B)$. From \eqref{lapF}, using $f_4=f+S^2z^2$, $f_3=S^{3/2}z$ and
$C=tS^{5/2}m$, we have
\[
 2A+B\le tS^3\bigl(\frac{f}{S^2}+(1-c)z^2+2czm+3df^*x\bigr).
\]
Also \eqref{ABlower} gives
\[
 -(5+4\delta/3)(A+2B)
 \le-tS^3(5+4\delta/3)
                  \left(3m^2+\frac{4f^*x^2}{x+2y}\right).
\]
The first term in \eqref{symprojection}, after division by $tS^3$, is
\[
 \frac{S^2}{3t(S^2+2f_4)}
        \left(3tm+4\delta\sum_i\mu_i^3b_i\right)^2.
\]
Substituting these estimates and $-4\delta f_3C/(tS^4)=-4\delta zm$
in \eqref{symprojection}, we obtain \eqref{Phi} after division by $tS^3$.

On the other hand, \eqref{eq:4.3} and Lemma~\ref{commutator} give
\begin{align*}
 |u|^2+\frac32Sf+(A+2B)
 &=t(2t-1)S^3+4(A-B)\\
 &\le t(2t-1)S^3+2tS^3=t(2t+1)S^3.
\end{align*}
Taking $\delta=0$ in \eqref{symprojection}, we have
$|u|^2\ge9S^3f^*y/(4-6y)$. Hence, by \eqref{ABlower},
\[
 \frac{3f}{2S^2}+3tm^2+
 f^*\left(\frac{4tx^2}{x+2y}+\frac{9y}{4-6y}\right)
 \le t(2t+1).
\]
This is equivalent to \eqref{Psi}.
\end{proof}

\begin{lemma}\label{monotone}
Let $0<\delta\le3/2$ and $0\le\gamma\le3/4$. With the limiting
quantities in \eqref{Phi} and \eqref{Psi} held fixed, we have, for $t>0$,
\[
 \frac{d}{dt}\bigl(\Phi_\delta(t)+\gamma\Psi(t)\bigr)\ge\frac14.
\]
In particular, if $0<t\le3/8$, then
$\Phi_\delta(t)+\gamma\Psi(t)
\le\Phi_\delta(3/8)+\gamma\Psi(3/8)$.
\end{lemma}
\begin{proof}
Taking $\vec w$ to be the $i$th coordinate vector in \eqref{eq:3.5}, we have
$\mathcal G_{ii}(1+2\mu_i^2)\le1$. If $\mu_i^2\le1/2$ for every $i$, then
\[
 \frac{A}{tS^3}=\sum_i\mu_i^2\mathcal G_{ii}
 \le\frac12\tr\mathcal G=\frac12.
\]
If $\mu_k^2>1/2$ for some $k$, then $\mathcal G_{kk}\le1/2$ and
$\mu_i^2\le1-\mu_k^2$ for $i\ne k$. Hence
\begin{align*}
 \frac{A}{tS^3}
 &\le\mu_k^2\mathcal G_{kk}
       +(1-\mu_k^2)(1-\mathcal G_{kk})\\
 &=1-\mu_k^2+(2\mu_k^2-1)\mathcal G_{kk}\le\frac12.
\end{align*}
Since $B\le A$, it follows from \eqref{ABlower} that
\[
m^2+\frac43f^*\frac{x^2}{x+2y}\leq\frac12.
\]
Differentiating \eqref{Phi} and \eqref{Psi} with respect to $t$ gives
\begin{align*}
 \frac{d}{dt}(\Phi_\delta+\gamma\Psi)
 ={}&\frac{(3-2\delta)(1+2\delta)f}{2t^2S^2}
 +\frac{16\delta^2S^2}{3t^2(S^2+2f_4)}
                         \left(\sum_i\mu_i^3b_i\right)^2\\
 &+\frac{8\delta^2}{3t^2}
                   \left(\frac{f}{2S^2}-\sum_i\mu_i^2b_i^2\right)
 +\frac{(3+2\delta)^2f^*y}{t^2(4-6y)}\\
 &+2-\frac{3S^2m^2}{S^2+2f_4}
 +\gamma\left(\frac23-2m^2
                   -\frac{8f^*x^2}{3(x+2y)}\right)
 +\frac{8\gamma t}{3}.
\end{align*}
By Lemma~\ref{spectral2},
$\frac{f}{2S^2}-\sum_i\mu_i^2b_i^2\ge0$ for $f>0$; this also holds
for $f=0$, since then $b_i=0$. As $0<\delta\le3/2$, the terms containing
$t^{-2}$ are nonnegative.
Since $(S^2+2f_4)/(3S^2)\ge1/3$ and $m^2\le1/2$, we have
$2-\frac{3S^2m^2}{S^2+2f_4}\ge1/2$. Moreover, the preceding bound on
$m^2+4f^*x^2/(3(x+2y))$ gives
\[
 \frac23-2m^2-\frac{8f^*x^2}{3(x+2y)}\ge-\frac13.
\]
Thus the derivative is at least
$1/2-\gamma/3\ge1/4>0$, proving the assertion.
\end{proof}

\section{Proof of Theorem~\ref{main}}\label{mainproof}
\begin{theorem}\label{gap}
Let $X:M^n\to\mathbb R^{n+1}$ be a complete self-shrinker with
constant squared norm $S$ of the second fundamental form. If $S\le8/5$,
then $S=0$ or $S=1$.
\end{theorem}
\begin{proof}
Suppose that $1<S\leq8/5$. By Proposition~\ref{newtensor} and Lemma~\ref{monotone}, we have
\[
0\leq\Phi_\delta(t)+\gamma\Psi(t)\leq\Phi_\delta(\frac38)+\gamma\Psi(\frac38)=:\Theta.
\]
We shall prove $\Theta<0$ in the two cases $0\le z^2\le2/5$ and $2/5\le z^2\le1$. Put
\[
Q_0=20+16\delta/3+\gamma,\quad B_1=Q_0-3dk_\delta,\quad
C_1=\frac83(3+2\delta)^2+6\gamma.
\]
For the values of $\delta$ and $\gamma$ chosen below, $0<B_1<Q_0$. Let
\[
\mathfrak r=\frac1{\sqrt6}(\sqrt{8Q_0}-\sqrt{8B_1}-\sqrt{C_1})_+.
\]
For $x,y\geq0$, $x+2y\leq1$, we have
\begin{align*}
3dk_\delta x-\frac{Q_0x^2}{x+2y}
&=2(Q_0+B_1)y-B_1(x+2y)-\frac{4Q_0y^2}{x+2y}\\
&\leq2y(\sqrt{Q_0}-\sqrt{B_1})^2.
\end{align*}
At $x=y=0$, both sides are zero. Put $a=\sqrt{Q_0}-\sqrt{B_1}$. Since
\[
 \frac{d}{dy}\left(2a^2y-\frac{C_1y}{4-6y}\right)
 =2a^2-\frac{4C_1}{(4-6y)^2},
\]
the function is decreasing on $[0,2/3)$ when $8a^2\le C_1$, and its maximum is zero.
If $8a^2>C_1$, the derivative vanishes at
$4-6y=\sqrt{2C_1}/a$; substitution gives the maximum
$(\sqrt8a-\sqrt{C_1})^2/6$. Thus the maximum is $\mathfrak r^2$. Since $0\le y\le1/2$, it follows that
\[
3dk_\delta x-\frac{Q_0x^2}{x+2y}-\frac{C_1y}{4-6y}\leq\mathfrak r^2.
\]
Therefore, adding \eqref{Phi} and $\gamma$ times \eqref{Psi}, we obtain
\begin{align}\label{endineq}
0\le\Theta\le{}&-\frac14+\frac7{16}\gamma
 +\left[\frac{16\delta}{3}(\delta-1/2)-\gamma\right]\frac{f}{S^2}\notag\\
 &+k_\delta(1-c)z^2+(2ck_\delta-4\delta)zm
 -(15+4\delta+3\gamma/4)m^2\notag\\
 &+\mathfrak r^2f^*
 -\frac{8S^2}{9(S^2+2f_4)}
        \left(\frac98m+4\delta\sum_i\mu_i^3b_i\right)^2\notag\\
 &-\frac{64\delta^2}{9}
             \left(\frac{f}{2S^2}-\sum_i\mu_i^2b_i^2\right).
\end{align}

\smallskip
\noindent\emph{Case 1: $0\leq z^2\leq2/5$.}
Take $\delta=1/15+2z^2/5$ and $\gamma=11/20-51z^2/50$.
For $f>0$, Lemma~\ref{spectral2} gives
\[
 \frac{f}{2S^2}-\sum_i\mu_i^2b_i^2
 \ge\frac{S^2}{S^2-f_4}\left(\sum_i\mu_i^3b_i\right)^2.
\]
Using this estimate and the identity
\[
 \frac{S^2+2f_4}{3S^2}
       +\frac{2(S^2-f_4)}{3S^2}=1,
\]
the Cauchy--Schwarz inequality yields
\begin{align*}
 &\frac{8S^2}{9(S^2+2f_4)}
           \left(\frac98m+4\delta\sum_i\mu_i^3b_i\right)^2
 +\frac{64\delta^2}{9}
           \left(\frac{f}{2S^2}-\sum_i\mu_i^2b_i^2\right)\\
 &\quad\ge
 \frac{8S^2\left(\frac98m+4\delta\sum_i\mu_i^3b_i\right)^2}
      {9(S^2+2f_4)}
 +\frac{64\delta^2S^2\left(\sum_i\mu_i^3b_i\right)^2}
      {9(S^2-f_4)}
 \ge \frac38m^2.
\end{align*}
Indeed, apply $u^2/a+v^2/b\ge(u-v)^2/(a+b)$ to the two terms,
whose numerators differ by $\frac98m$ after putting them over denominators
$\frac98(S^2+2f_4)/S^2$ and $\frac94(S^2-f_4)/S^2$.
The sum of these denominators is $\frac{27}{8}$.
If $f=0$, then $b_i=0$ for every $i$ and $f_4\le S^2$, so the same
conclusion follows directly from the first term.
The remaining terms containing $m$ satisfy
\begin{align*}
 &(2ck_\delta-4\delta)zm
 -(15+4\delta+3\gamma/4+\frac38)m^2\\
 &=-\left(15+4\delta+3\gamma/4+\frac38\right)
 \left(m-\frac{(ck_\delta-2\delta)z}
                 {15+4\delta+3\gamma/4+\frac38}\right)^2
 +\frac{(ck_\delta-2\delta)^2z^2}{15+4\delta+3\gamma/4+\frac38}.
\end{align*}
Discarding this nonpositive square in \eqref{endineq}, we obtain
\begin{align}\label{caseonebound}
\Theta\le{}&-\frac14+\frac7{16}\gamma
 +z^2\left[k_\delta(1-c)
 +\frac{(ck_\delta-2\delta)^2}{15+4\delta+3\gamma/4+3/8}\right]\notag\\
 &-\left[\gamma-\frac{16\delta}{3}(\delta-1/2)\right]\frac{f}{S^2}
 +\mathfrak r^2f^*.
\end{align}
We estimate the expression in square brackets in the first line directly.
Substituting $\delta=1/15+2z^2/5$, $\gamma=11/20-51z^2/50$
and $c=137/100-z^2/10$ gives
\begin{align*}
 &\left(15+4\delta+\frac{3\gamma}{4}+\frac38\right)
 \left[\frac{401}{1000}+\frac{21z^2}{200}-k_\delta(1-c)
       -\frac{(ck_\delta-2\delta)^2}{15+4\delta+3\gamma/4+3/8}\right]\\
 &=\frac{536287+35179589z^2-82750538z^4
                    -6011520z^6-1843200z^8}{162000000}.
\end{align*}
Since $0\le z^2\le2/5$, we have $z^6\le(2/5)z^4$ and
$z^8\le(4/25)z^4$. Hence the numerator is at least
\[
536287+z^2(35179589-85450058z^2)
\ge536287+\frac{4997829}{5}z^2>0.
\]
The factor $15+4\delta+3\gamma/4+3/8$ is positive. Therefore
\[
k_\delta(1-c)+\frac{(ck_\delta-2\delta)^2}{15+4\delta+3\gamma/4+3/8}
\le\frac{401}{1000}+\frac{21z^2}{200}.
\]
It follows that the first line of \eqref{caseonebound} is at most
\begin{align*}
-\frac14+\frac7{16}\gamma
 +z^2\left(\frac{401}{1000}+\frac{21z^2}{200}\right)
 &=-\frac3{320}-\frac{181}{4000}z^2+\frac{21}{200}z^4\\
 &\le-\frac3{320}-\frac{13}{4000}z^2<0.
\end{align*}
Here we used $z^4\le(2/5)z^2$.
It remains to show that
\[
\mathfrak r^2<\gamma-\frac{16\delta}{3}(\delta-1/2)
=\frac{1901}{2700}-\frac{107z^2}{450}-\frac{64z^4}{75}.
\]

Next, we have
\begin{align*}
Q_0&=3763/180+167z^2/150,\\
B_1&=16807/4500+367z^2/750+24z^{4}/25,\\
C_1&=39799/1350+1631z^2/225+128z^{4}/75.
\end{align*}
In particular, $Q_0>B_1>0$, $\gamma-16\delta(\delta-1/2)/3>0$. Directly from these polynomials, $8B_1,C_1\in[(27/5)^2,(29/5)^2]$ and $6\left(\frac{1901}{2700}-\frac{107z^2}{450}-\frac{64z^4}{75}\right)\in[(8/5)^2,(21/10)^2]$. For $a^2\le u\le b^2$ with $0<a<b$, concavity of the square root gives
\[
 \sqrt u\ge\frac{u+ab}{a+b}.
\]
Also $\sqrt u\le(u+13^2)/26$, by $(\sqrt u-13)^2\ge0$.
Apply the first inequality with $(a,b)=(27/5,29/5)$ to both $8B_1$ and
$C_1$, and with $(a,b)=(8/5,21/10)$ to $6\left(\frac{1901}{2700}-\frac{107z^2}{450}-\frac{64z^4}{75}\right)$. Apply the second to $8Q_0$.
We obtain
\begin{align*}
&\sqrt{8B_1}+\sqrt{C_1}+\sqrt{6\left(\frac{1901}{2700}-\frac{107z^2}{450}-\frac{64z^4}{75}\right)}-\sqrt{8Q_0}\\
&\quad\geq\frac{5}{56}\left(8B_1+C_1+\frac{1566}{25}\right)\\
&\qquad+\frac{10}{37}\left(6\left(\frac{1901}{2700}-\frac{107z^2}{450}-\frac{64z^4}{75}\right)+\frac{84}{25}\right)-\frac{4Q_0}{13}-\frac{13}{2}\\
&\quad=\frac{-19843200z^{4}+9767130z^2+371459}{36363600}>0.
\end{align*}
The numerator of the last expression is a concave quadratic in $z^2$, with endpoint values $371459$ and $1103399$ at $z^2=0$ and $z^2=2/5$, respectively. It is therefore positive throughout the interval. Consequently $\mathfrak r^2<\gamma-16\delta(\delta-1/2)/3$.
Since $0\le f^*\le f/S^2$, the second line of \eqref{caseonebound}
is nonpositive. Together with the estimate of its first line, this gives
\[
0\le\Theta\le-\frac3{320}-\frac{13}{4000}z^2<0,
\]
which is impossible.

\smallskip
\noindent\emph{Case 2: $2/5\leq z^2\leq1$.}
Take $\delta=3z^2/2$, $\gamma=0$, and put $\omega=3/2-5\delta/6\in[1/4,1]$. We have
\[
Q_0=20+16\delta/3,\qquad B_1=89/25-242\delta/75+8\delta^2/5.
\]
In particular $0<B_1<Q_0$. On $[3/5,3/2]$,
\[
\sqrt{Q_0}\leq4.82+\frac59(\delta-3/5),\qquad
\sqrt{B_1}\geq1.67-\delta/3>0.
\]
Indeed, the quantities on both sides are positive, and
\begin{align*}
 \left(4.82+\frac59(\delta-3/5)\right)^2-Q_0
 &=\frac{25}{81}\left(\delta-\frac{141}{250}\right)^2+\frac4{125}>0,\\
 B_1-(1.67-\delta/3)^2
 &=\frac{67}{45}\left(\delta-\frac{951}{1340}\right)^2
       +\frac{887}{41875}>0.
\end{align*}
Consequently,
\begin{align*}
 \sqrt{8Q_0}-\sqrt{8B_1}-\sqrt{C_1}
 &\le\sqrt8\left(\frac{169}{60}+\frac{8\delta}{9}\right)
       -\sqrt{\frac83}(3+2\delta)\\
 &\le\frac{37}{12}-\frac{13\delta}{18}.
\end{align*}
The last inequality follows by evaluating the two affine functions at
$\delta=3/5$ and $\delta=3/2$. Since the last expression is positive
on this interval, the definition of $\mathfrak r$ gives
\begin{equation}\label{rbound}
 \frac{\mathfrak r}{\sqrt\delta}
 \le\frac{37/12-13\delta/18}{\sqrt{6\delta}}.
\end{equation}

Expanding
\[
 \left(\frac38m+\frac{4\delta}{3}\sum_i\mu_i^3b_i
       -\frac{\omega z(S^2+2f_4)}{8S^2}\right)^2
\]
and using its nonnegativity, we obtain
\begin{align*}
 &-\frac{8S^2}{9(S^2+2f_4)}
            \left(\frac98m+4\delta\sum_i\mu_i^3b_i\right)^2\\
 &\quad\le-2\omega z\left(\frac38m+\frac{4\delta}{3}\sum_i\mu_i^3b_i\right)
              +\frac{\omega^2z^2(S^2+2f_4)}{8S^2}.
\end{align*}
The resulting $m$-terms satisfy
\begin{align*}
 &2(ck_\delta-2\delta-\frac38\omega)zm-(15+4\delta)m^2\\
 &\le\frac{(ck_\delta-2\delta-\frac38\omega)^2z^2}{15+4\delta}.
\end{align*}
For $f>0$, use \eqref{projectbound} with $a_0=\mathfrak r$ and
$b_0=8\delta/3$ to estimate the remaining terms by
\[
 \mathfrak r^2f^*-\frac{64\delta^2}{9}(\frac{f}{2S^2}-\sum_i\mu_i^2b_i^2)
 \le \frac{f}{S^2}\left(\mathfrak r-\frac{8\delta}{3}|r-cz|\right)_+^2
       -\frac{64\delta^2}{9}\frac{f}{S^2}(\frac{f_4}{S^2}-\frac12).
\]
In addition,
\[
 \frac38\omega^2z^2\left(\frac{S^2+2f_4}{3S^2}\right)
 =\frac{\omega^2z^2(1+2z^2)}8+\frac{\omega^2z^2}{4}\frac{f}{S^2},
 \qquad -\frac{8\delta\omega}{3}z\left(\sum_i\mu_i^3b_i\right)
 =-\frac{8\delta\omega}{3}zr\frac{f}{S^2}.
\]
Using $z^2=2\delta/3$, the coefficient of $f/S^2$ becomes
\begin{align*}
 &\frac{16\delta}{3}(\delta-1/2)
   -\frac{64\delta^2}{9}(z^2-1/2)+\frac{\omega^2z^2}{4}\\
 &=\delta\left((1-z^2)(64\delta/9-8/3)+\omega^2/6\right).
\end{align*}
Finally apply \eqref{rbound} to the positive-part term. We obtain
\begin{equation}\label{Vineq}
\begin{aligned}
\Theta\le{}&-\frac14+z^2\left[k_\delta(1-c)
 +\frac{(ck_\delta-2\delta-3\omega/8)^2}{15+4\delta}
 +\frac{\omega^2(1+2z^2)}8\right]\\
 &+\delta\frac{f}{S^2}\left[V(r)-\frac{64\delta}{9}\frac{f}{S^2}\right],
\end{aligned}
\end{equation}
where $r=\frac{S^2}{f}\sum_i\mu_i^3b_i$ if $f>0$ and
\begin{align*}
V(r)={}&(1-z^2)(64\delta/9-8/3)+\omega^2/6-8\omega zr/3\\
&+\left(\frac{37/12-13\delta/18}{\sqrt{6\delta}}
             -\frac{8\sqrt\delta}{3}|r-cz|\right)_+^2.
\end{align*}
For $f=0$, the same argument gives \eqref{Vineq} with its last line omitted.

For $f>0$, $r=z+\frac{S^2}{f}\sum_i\mu_i b_i^2$. If $\mu_i<0$, then $z\leq(1-\mu_i^2)^{3/2}$, and consequently $z^2\mu_i^2\leq27/256<1/9$. Hence $r\geq z-1/(3z)$. Indeed, writing $u=\mu_i^2\in[0,1]$ for a negative $\mu_i$ gives
$z^2u\le u(1-u)^3\le27/256<1/9$; the maximum occurs at $u=1/4$.
Thus every $\mu_i\ge-1/(3z)$, and $r=z+\frac{S^2}{f}\sum_i\mu_i b_i^2\ge z-1/(3z)$.
On $[z-1/(3z),cz]$, the positive-part term in $V$ is the square of the positive
part of an increasing affine function, and is convex. The remaining
$r$-term is affine. For $r\ge cz$, both the square term and $-8\omega zr/3$ are
nonincreasing. It therefore suffices to prove, at $r=z-1/(3z)$ and $r=cz$,
\begin{equation}\label{Vbound}
V(r)<(52\omega-1)/30.
\end{equation}
Since $(c-1)z^2\geq0.132$, we have
\[
\frac{64\delta}{9}\left((c-1)z+\frac1{3z}\right)^2>\frac{32}{3}(0.13+1/3)^2>2>
\frac{(37/12-13\delta/18)^2}{6\delta}.
\]
The last inequality follows because the expression on its right is decreasing
for $3/5\le\delta\le3/2$ and has value $2809/1440<2$ at $\delta=3/5$. Thus the positive part in $V(z-1/(3z))$ vanishes. Substitution gives
\[
(52\omega-1)/30-V(z-1/(3z))=\frac{2037\omega^2-2596\omega+897}{450}>0,
\]
by the negative discriminant. At $r=cz$, direct expansion yields
\begin{align*}
6\delta\big[(52\omega-1)/30-V(cz)\big]
=\frac1{90000}\big(&110592\omega^4-2919456\omega^3\\
&+6706592\omega^2-5535640\omega+1651295\big).
\end{align*}
The expression on the right can also be written as
\begin{align*}
\frac{256}{81}\bigg[&\frac{5701}{800}(1-\omega)^4
 +\frac{23259}{4000}(\omega-1/4)(1-\omega)^3\\
 &+\frac{7943}{2500}(\omega-1/4)^2(1-\omega)^2
 +\frac{21243}{5000}(\omega-1/4)^3(1-\omega)\\
 &+\frac{1487}{10000}(\omega-1/4)^4\bigg].
\end{align*}
All five terms are nonnegative for $1/4\le\omega\le1$.
They cannot all vanish: the first is positive at $\omega=1/4$,
the last is positive at $\omega=1$, and both are positive in the interior.
This proves \eqref{Vbound} at $r=cz$. The preceding convexity and
monotonicity then give $V(r)\le(52\omega-1)/30$ for all $r\ge z-1/(3z)$.

By \eqref{Vineq} and \eqref{Vbound},
\begin{align*}
 \delta \frac{f}{S^2}\left[V(r)-\frac{64\delta}{9}\frac{f}{S^2}\right]
 &\le\frac{52\omega-1}{30}\delta \frac{f}{S^2}-\frac{64}{9}(\delta \frac{f}{S^2})^2\\
 &=-\frac{64}{9}\left(\delta \frac{f}{S^2}-\frac{3(52\omega-1)}{1280}\right)^2
      +\frac{(52\omega-1)^2}{25600}.
\end{align*}
Discarding the square and substituting in \eqref{Vineq}, we obtain
\begin{align}\label{casetwobound}
\Theta\le{}&-\frac14+z^2\bigg[
 \frac{100+12\omega+32\omega^2}{375}
 +\frac{\omega^2(329-256\omega)^2}{10^6(15+4\delta)}\notag\\
 &\hspace{25mm}+\frac{3\omega}{8}
       \left(\frac\omega3(1+4\delta/3)-\frac23\right)\bigg]
 +\frac{(52\omega-1)^2}{25600}.
\end{align}
Here we used $\delta=9/5-6\omega/5$, $z^2=6/5-4\omega/5$,
$c=137/100-z^2/10$ and $k_\delta=4+8\delta/3$ to expand
\[
 k_\delta(1-c)+\frac{(ck_\delta-2\delta-3\omega/8)^2}{15+4\delta}
 +\frac{\omega^2(1+2z^2)}8.
\]
The same bound holds when $f=0$.
The remaining rational term can be bounded without introducing another
auxiliary function. Completing the square gives
\[
 \omega(329-256\omega)
 =\frac{329^2}{1024}-256\left(\omega-\frac{329}{512}\right)^2<106.
\]
Since $1/4\le\omega\le1$, its left side is positive. Also
$15+4\delta\ge87/5$, because $\delta\ge3/5$. Therefore
\[
 0\le\frac{\omega^2(329-256\omega)^2}{10^6(15+4\delta)}
 <\frac{106^2}{10^6(87/5)}<\frac1{1000}.
\]
Substitute the preceding bound, $z^2=6/5-4\omega/5$ and
$\delta=9/5-6\omega/5$ in \eqref{casetwobound}.
Collecting powers of $1-\omega$ gives
\begin{align*}
\Theta\le{}&\frac4{25}(1-\omega)^4
 +\frac{31}{3750}(1-\omega)^3-\frac{739}{8000}(1-\omega)^2\\
 &-\frac{61}{240000}(1-\omega)-\frac{563}{128000}\\
\le{}&\frac{153}{40000}(1-\omega)^2-\frac{563}{128000}
\le-\frac{719}{320000}<0.
\end{align*}
Here we used $0\le1-\omega\le3/4$ and
\[
 \frac4{25}(1-\omega)^4\le\frac9{100}(1-\omega)^2,
 \qquad
 \frac{31}{3750}(1-\omega)^3\le\frac{31}{5000}(1-\omega)^2,
\]
together with $-61(1-\omega)/240000\le0$ and
$9/100+31/5000-739/8000=153/40000$. Moreover,
\[
 \frac{153}{40000}(1-\omega)^2-\frac{563}{128000}
 \le\frac{1377}{640000}-\frac{2815}{640000}
 =-\frac{719}{320000}<0.
\]
This contradicts $\Theta\ge0$ and completes Case 2.
Thus $1<S\le8/5$ is impossible. For $S\le1$, \eqref{eq:2.1} gives
$0\le\sum h_{ijk}^2=S(S-1)\le0$, so $S=0$ or $S=1$.
\end{proof}

\begin{proof}[Proof of Theorem~\ref{main}]
Suppose $S>8/5$. Since
$R=H^2-S\ge-8/5$, we have
\[
 0<S-\frac85\le H^2\le nS.
\]
The mean curvature vector is therefore nowhere zero. Choosing the unit
normal in its direction, we have $H>0$.
The function $S/H^2$ is bounded above by $S/(S-8/5)$.
The Gauss equation also gives, in a principal frame,
\[
 \operatorname{Ric}(e_i,e_i)=H\lambda_i-\lambda_i^2
 \ge-(\sqrt n+1)S.
\]
We may therefore apply the generalized maximum principle of Cheng and
Peng \cite{CP} to $S/H^2$.

The quotient identity in \cite[(2.9)]{CLW1} reads
\begin{equation}\label{eq:3.23}
 \frac12\mathcal L\frac{S}{H^2}
 =\frac1{H^4}\sum_{ijk}(h_{ij}\nabla_kH-Hh_{ijk})^2
 -\frac1H\left\langle\nabla H,\nabla\frac{S}{H^2}\right\rangle.
\end{equation}
Choose $p_m\in M$ such that
\[
 \frac{S}{H^2}(p_m)\longrightarrow\sup_M\frac{S}{H^2},\qquad
 \left|\nabla\frac{S}{H^2}\right|(p_m)\longrightarrow0,
 \qquad\limsup_m\mathcal L\frac{S}{H^2}(p_m)\le0.
\]
Since $S$ is constant,
\[
 \nabla\frac{S}{H^2}=-\frac{2S}{H^3}\nabla H,
 \qquad
 |\nabla H|\le\frac{(nS)^{3/2}}{2S}
                       \left|\nabla\frac{S}{H^2}\right|.
\]
Consequently
\begin{equation}\label{eq:3.24}
 |\nabla H|(p_m)\longrightarrow0.
\end{equation}
Because $H\ge\sqrt{S-8/5}>0$, the last term in \eqref{eq:3.23}
tends to zero. The first term is nonnegative, so the maximum principle implies
\[
 \frac1{H^4}\sum_{ijk}(h_{ij}\nabla_kH-Hh_{ijk})^2(p_m)
 \longrightarrow0.
\]
The upper bound $H^2\le nS$ then implies that its numerator tends to zero.
Hence
\begin{align*}
 H^2\sum_{ijk}h_{ijk}^2
 &\le2\sum_{ijk}(Hh_{ijk}-h_{ij}\nabla_kH)^2
       +2S|\nabla H|^2\longrightarrow0.
\end{align*}
Since $H^2\ge S-8/5>0$, we obtain
$\sum_{ijk}h_{ijk}^2(p_m)\to0$, whereas \eqref{eq:2.1} gives
\[
 \sum_{ijk}h_{ijk}^2=S(S-1)>0
\]
on $M$. This contradiction proves $S\le8/5$.

By Theorem~\ref{gap}, $S=0$ or $S=1$. Equation~\eqref{eq:2.1} then
implies $h_{ijk}=0$. Thus the second fundamental form is parallel.
By the classification in \cite{l} and the completeness of $M$, the image
$X(M)$ is a hyperplane, a sphere, or a round cylinder.
The self-shrinker equation places the hyperplane through the origin and
centers the spherical factor at the origin in its normal Euclidean subspace.
For a spherical factor of dimension $k$ and radius $r$, that equation gives
$r^2=k$. Thus, up to an orthogonal transformation,  the image $X(M^n)$
is one of  $S^{n}(\sqrt{n})$,  $\mathbb{R}^{n}$, 
$S^k (\sqrt k)\times \mathbb{R}^{n-k}$, $1\leq k\leq n-1$.
We finish our proof of Theorem \ref{main}.
\end{proof}

 \vskip2mm
\noindent
{\bf Acknowledgement}.
This work was partly supported by MEXT Promotion of Distinctive Joint Research Center Program JPMXP0723833165 
and Osaka Metropolitan University Strategic Research Promotion Project (Development of International Research Hubs).

\end{document}